\documentclass{amsart}

\usepackage{amsmath,amsfonts,amssymb,amsthm,amscd,latexsym,euscript}

\usepackage[all]{xy}

\newcommand{\op}{{\ensuremath{\textup{op}}}}

\newcommand{\dgrm}[1]{\ensuremath{\smash{\underset{\widetilde{\hphantom{#1}}}{#1}} \mathstrut}}

\newcommand{\cal}[1]{\ensuremath{\mathcal #1}}

\newtheorem {theorem1}{Theorem}[section]
\newtheorem {theorem}[theorem1]{Theorem}
\newtheorem {corollary}[theorem1]{Corollary}
\newtheorem {proposition}[theorem1]{Proposition}
\newtheorem {lemma}[theorem1]{Lemma}
\theoremstyle{definition}
\newtheorem {definition}[theorem1]{Definition}

\theoremstyle{remark}
\newtheorem {remark}[theorem1]{Remark}

\newcommand{\cat}[1]{\ensuremath{\EuScript #1}}

\DeclareMathOperator{\map}{\textup{map}}

\DeclareMathOperator{\Spec}{\ensuremath{\textup{Sp}}}

\newcommand{\holim}{\ensuremath{\mathop{\textup{holim}}}}
\newcommand{\colim}{\ensuremath{\mathop{\textup{colim}}}}
\newcommand{\hocolim}{\ensuremath{\mathop{\textup{hocolim}}}}

\def\id{\ensuremath{{\rm id}}}

\newcommand{\rarrow}{\rightarrow}
\newcommand{\Id}{\ensuremath{\textup{Id}}}

\newcommand{\overcat}{\ensuremath{\!\downarrow\!}}

\newcommand{\Rosicky}{Rosick\'y }

\newcommand{\co}{\colon\thinspace}
\newcommand{\sm}[2]{\ensuremath{#1\smash[b]{\wedge\,}#2}}
\newcommand{\aus}{\raisebox{1pt}{\ensuremath{\,{\scriptstyle\in}\,}}}%

\newcommand{\ul}[1]{\underline{#1}}

\newcommand{\diag}[2]{ \begin{align} \begin{split} \xymatrix{#1} \end{split} \label{#2} \end{align}}%

\newcommand{\diagr}[1]{ \begin{equation*} \xymatrix{#1} \end{equation*}}%

\newcounter{zahl}%
\newenvironment{punkt}{\begin{list}{{\rm{(\roman{zahl})}}}%
    {\usecounter{zahl}%
     \setlength{\leftmargin}{0pt} \setlength{\itemindent}{4pt} \setlength{\topsep}{2pt} \setlength{\parsep}{2pt} }}%
    {\end{list}}%

\usepackage{hyperref} 

\begin{document}

\SelectTips{cm}{10}

\title  [Calculus of functors]
        {Calculus of functors and model categories}
\author{Georg Biedermann}
\author{Boris Chorny}
\author{Oliver R\"ondigs}

\thanks{The work was supported by the SFB 701 at the University of Bielefeld.}

\address{Department of Mathematics, Middlesex College, The University of Western Ontario, London, Ontario N6A 5B7, Canada}
\address{D-MATH, ETH Zentrum, 8092 Z\"urich, Switzerland}
\address{Fakult\"at f\"ur Mathematik, Universit\"at Bielefeld, Postfach 100 131, D-33501 Bielefeld, Germany}

\email{gbiederm@uwo.ca}
\email{chorny@math.ethz.ch}
\email{oroendig@math.uni-bielefeld.de}



\date{\today}
\dedicatory{}
\commby{}

\begin{abstract}
The category of small covariant functors from simplicial sets to 
simplicial sets supports the projective model structure \cite{Chorny-Dwyer}. 
In this paper we construct various localizations of the projective model 
structure and also give a variant for functors from simplicial sets to 
spectra. We apply these model categories in the study of calculus of 
functors, namely for a classification of polynomial and homogeneous functors.  
In the $n$-homogeneous model structure, the
$n$-th derivative is a Quillen
functor to the category of spectra with $\Sigma_n$-action. After taking into account only 
finitary functors -- which may be done in two different ways -- 
the above Quillen map becomes a Quillen equivalence. 
This improves the classification of finitary homogeneous functors by
T.~G.~Goodwillie \cite{Goo:calc3}.
\end{abstract}

\maketitle

\section{Introduction}
Calculus of homotopy functors applies to functors from spaces to spaces 
or spectra which preserve weak equivalences. It may be viewed as an interpolation between 
stable and unstable homotopy theory by analyzing carefully the rate of 
change of such functors. Developed around 1990 by
Thomas G.~Goodwillie, calculus of functors has had spectacular applications to
geometric topology \cite{Goo:calc1, Goo:calc2} and homotopy theory
\cite{Arone-Mahowald}. Although at the present time calculus of
functors is a well developed and ramified theory, foundations of the
subject remain technically involved.

In the current work we introduce a categorical approach to these
foundations. In order to overcome set-theoretical difficulties,
we consider only functors from spaces to spaces or from spaces
to spectra which are determined by their restriction to some 
small subcategory. 
We suggest to implement the
machinery developed by Goodwillie as a part of a simplicial model category
structure on this functor category. 
As an immediate advantage of this approach we obtain well-behaved
mapping spaces between functors.
For technical reasons, we use simplicial sets instead of topological spaces.
This is justified by Kuhn's overview 
article \cite{Kuhn:overview}, 
where first steps to an axiomatization of the theory are taken.
Finally, all functors are assumed to be simplicial (or continuous, or enriched).

The projective model structure, in which weak equivalences and fibrations
are detected objectwise, was constructed in
\cite{Chorny-Dwyer}. In this paper we present several new model
structures on the category of small functors, and each of these
reflects a certain aspect of Goodwillie's calculus. 

After necessary preliminaries on small functors in
Section~\ref{prelim} we construct in Section~\ref{homotopy} a
localization of the projective model structure such that the new
fibrant objects are precisely the objectwise fibrant homotopy
functors. This is the starting point for calculus of functors,
since Goodwillie's machinery is intended for homotopy functors only. 

In Section~\ref{n-excisive structure}, we localize the homotopy
model structure on the category of small functors from spaces to
spaces. The new fibrant objects are precisely the $n$-excisive
fibrant homotopy functors. This result may be viewed as a
classification of $n$-polynomial functors. Goodwillie's $n$-th
polynomial approximation $P_n$
is equivalent to a fibrant replacement in our $n$-excisive model
structure. An immediate advantage of having a model category structure 
is that the fibrant replacement (equivalent to $P_n$) is universal up 
to homotopy with respect to maps into an arbitrary $n$-excisive functor. 
This is an improvement of Goodwillie's result, which verifies 
the universal property only on the level of the homotopy category \cite[1.8]{Goo:calc3}.

In the simpler category of functors from finite pointed
spaces to all pointed spaces, Lydakis has constructed the homotopy model
structure as well as the 1-excisive (or stable) model 
structure (see \cite{Lydakis}, as well as its 
generalization \cite{DRO:enriched} to more general model categories).
Our work may be seen as a two-fold generalization of this work, 
since our results immediately apply to Lydakis' category. However,
there are plenty of interesting small functors which are not
determined by their values on finite spaces -- for example,
non-smashing Bousfield localizations. 
 
Another predecessor of our $n$-excisive model structure on the category of small functors is the $n$-excisive model structure constructed by W.~G.~Dwyer \cite{Dwyer_localizations} on the category of functors from finite $CW$-complexes to topological spaces.

In Section~\ref{section:spectrum-valued} we establish the stable 
projective, stable homotopy, and stable $n$-excisive model 
structures for small functors from (pointed) spaces to spectra. 
Then we recall and adapt several important definitions in
Section~\ref{section:Taylor-tower}.
In Section~\ref{section:homogeneous} we colocalize the
stable $n$-excisive model structure in order to obtain the
$n$-homogeneous model structure. In this model structure, the
fibrant and cofibrant objects are precisely those projectively
fibrant and cofibrant homotopy functors which are $n$-homogeneous. This model
structure may also be considered as a way to classify the
$n$-homogeneous functors up to homotopy. T.~Goodwillie has found
another, simpler classification, but it applies only for finitary $n$-homogeneous functors or for a restriction of an arbitrary functor to finite spaces. 
Any such functor is determined by its $n$-th derivative,
which is a spectrum with $\Sigma_n$-action. The construction of the derivative, as well 
as its interpretation as a Quillen functor, may be found in Section~\ref{section:F-equiv}. 


In the final Section~\ref{section:finitary}, we strengthen Goodwillie's
classification by introducing a finitary version of our
$n$-homo\-ge\-neous model structure and an $n$-homo\-ge\-neous model structure on the
category of functors from pointed finite simplicial sets to spectra. We prove that
the $n$-th derivative is a Quillen equivalence between this model category and the projective model
structure on the category of spectra with $\Sigma_n$-action. \\

Let \cal{S} denote the category of simplicial sets, and $\cal{S}_*$ the category
of pointed simplicial sets. The category of Bousfield-Friedlander spectra is denoted
$\Spec$. We stick to the common inconsistency in calling a category \cal{C} enriched in simplicial sets 
simply a simplicial category. The simplicial set of morphisms from $A$ to $B$ in a simplicial category
$\mathcal{C}$ is denoted by $\map_\mathcal{C}(A,B)$. If \cal{C} is cotensored over \cal{S}, the cotensor is denoted
by $A^K$ for $K\in\cal{S}$ and $A\in \mathcal{C}$. \\

\subsection*{Acknowledgments} We thank Rick Jardine for many helpful conversations on the 
subject matter of this paper. We thank also Andrew Mauer-Oats and the referee for numerous suggestions that helped to improve the exposition.

\section{Preliminaries on small functors}\label{prelim}
The object of study of this paper is the homotopy theory of simplicial functors from 
simplicial sets to simplicial sets or spectra. There are several cases including pointed 
and unpointed versions. We will first focus on endofunctors of unpointed simplicial sets.
The totality of these functors does not form a category in the usual sense -- natural 
transformations between two functors need not form a set in general.
Instead, we introduce a collection of functors which is on the one hand large enough to contain plenty of interesting functors,
and on the other hand small enough to form a category in the usual sense with small morphism sets. 

\begin{definition}
Let \cat K be a simplicial category. Any object $A\in \cat{K}$
defines a \emph{representable functor\/} 
\[ R^A\colon \cat K \rightarrow \cal S, \quad B \mapsto \map_{\cat{K}}(A,B).\]
A functor
$\dgrm X\colon \cat K \rarrow \cal S$ is called \emph{small} if \dgrm
X is a small weighted colimit of representable functors. 
We denote the category of small functors as $\cal S^\cat K$. 
\end{definition}

Since any representable functor is simplicial in the sense that it comes with functorial maps
\[ \map_{\cat{K}}(B,C) \rightarrow \map_{\cal{S}}\bigl(\map_{\cat{K}}(A,B),\map_{\cat{K}}(A,C)\bigr),\]
any small functor is simplicial or ``enriched over simplicial sets''.
Consequently all colimits and left Kan extensions in this article are to be taken in the enriched sense.
See Kelly's book  \cite{Kelly} 
for the necessary background on enriched category theory. A useful characterization of small functors is 
proved in \cite[Prop. 4.83]{Kelly}: 
a functor is small if and only if it is a left Kan extension from its restriction to a full small subcategory.
Note that Kelly calls small 
functors \emph{accessible} and weighted colimits \emph{indexed}. 

In other words, a small functor is determined by its values on some full small subcategory. 
Kelly proves also that small functors form a category enriched in
simplicial sets, where the simplicial mapping spaces are computed
using the formula \cite[4.41]{Kelly}. 
This category is closed under small weighted colimits \cite[Prop.~5.34]{Kelly}. That allows us to talk about mapping 
spaces $\map_{\cal S^\cat K}(\dgrm X, \dgrm Y)$ for all $\dgrm{X},\dgrm{Y} \in \cal S^\cat K$.  
The existence of weighted colimits implies in particular that $\cal S^\cat K$ is tensored over \cal S, as 
the functor $-\otimes K$ is the colimit over the trivial category weighted by $K\in\cal S$. 

Thus small functors provide a solution to the problem of writing down an honest category of functors.
For homotopy-theoretic constructions, the existence of certain limits is required. 
It turns out that, under some conditions on \cat{K}, the category of small functors $\cal S^\cat K$ is complete. 
The history of the problem is long, and work of Freyd \cite{Freyd} and \Rosicky \cite[Lemma 1]{Rosicky} provides a 
full answer for the question when the
category of small set-valued functors from a large category is complete. 
The work \cite{Lack} partly generalizes the results to the enriched settings, 
by showing that the category of small functors from \cat{K} to \cal S is complete if \cat{K} is a cocomplete simplicial category. 
The existence of weighted limits implies that $\cal S^{\cat K}$ is cotensored over \cal S, 
as the functor $(-)^K$ may be viewed as a limit over the trivial
category weighted by $K\in \cal S$.
By \cite[3.8]{Kelly} the functor $(-)^K$ is the right adjoint to the
tensor functor $-\otimes K$.
The results of B.~Day and S.~Lack allow us to consider simplicial model structures on the category of small functors $\cal S^{\cat K}$. 
The simplest model structure is the projective model structure
established in \cite{Chorny-Dwyer}.

\begin{definition}\label{projective structure}
A morphism $\dgrm{X}\to\dgrm{Y}$ in the category of small functors $\cal S^{\cat K}$ is 
\begin{punkt}
    \item
an \emph{objectwise equivalence} if $\dgrm{X}(K)\to\dgrm{Y}(K)$ is a weak equivalence in \cal{S} for all $K\in\cat{K}$.
    \item
an \emph{objectwise fibration} if $\dgrm{X}(K)\to\dgrm{Y}(K)$ is a fibration in \cal{S} for all $K\in\cat{K}$.
\end{punkt}
\emph{Projective cofibrations} are defined by the left lifting property with respect to trivial objectwise fibrations.
These classes form the \emph{projective model structure} on~$\cal{S}^{\cat{K}}$.
\end{definition}

The  (trivial) fibrations in the projective model structure are detected by mapping out of the enriched representable functors 
$R^A$ using the enriched Yoneda lemma 
\[ \map_{\cal S^{\cat K}}(R^A,X) \cong X(A). \]
Hence we obtain the following classes of generating (trivial) cofibrations:
\begin{align}\label{IandJ}
I &= \{R^A\otimes \partial \Delta^n \hookrightarrow R^A\otimes \Delta^n | A\in \cat K,\, n\geq 0\} \\
J &= \{R^A\otimes \Lambda^n_k\, \tilde\hookrightarrow R^A\otimes \Delta^n | A\in \cat K,\, n>0,\,0\leq k\leq n\}.
\end{align}
These are proper classes as soon as $\cat K$ is not a small category.
In this case, the factorization axioms are proved by a generalized small object argument \cite{pro-spaces}, 
which accepts certain classes of generating (trivial) cofibrations as input. 
The extra condition these classes have to satisfy 
is the co-solution set condition -- see \cite[3.1]{Chorny-Dwyer}.

\begin{definition}
Following \cite{pro-spaces} we call a model category \emph{class-cofibrantly generated} if there are 
two classes $I$ and $J$ of morphisms 
which admit the generalized small object argument and generate the model structure in the usual sense: 
$I\text{-inj}=\{\text{trivial fibrations}\}$ and $J\text{-inj}=\{\text{fibrations}\}$. 
\end{definition}

\begin{theorem}
The projective model category structure on $\cal S^{\cat K}$ is
simplicial, proper and class-cofibrantly generated.
\end{theorem}




Since the co-solution set for the small 
object argument \emph{is not chosen functorially} \cite[3.1]{Chorny-Dwyer}, 
the factorizations provided by the generalized small object argument are not functorial. 
So we traded functorial factorization for the possibility to work with a the larger class of all small functors. 
The reader who prefers 
to work with functors defined on a fixed small category --  say finite simplicial sets -- may do so. 
In this situation, the projective model structure and its localization to 
the $n$-excisive structure (to be described in Section~\ref{n-excisive structure}) are cofibrantly generated and 
therefore have functorial factorization. This was already known to Manos Lydakis (see \cite[p.~2]{Lydakis}).

The category of small functors has another important property: it is closed under composition. We will need this property in the next section. 

\begin{lemma} \label{composition is small}
The category of small functors $\cal S^{\cal S}$ is closed under composition. 
\end{lemma}

\begin{proof}
Given two small functors $\dgrm X, \dgrm Y\in \cal S^{\cal S}$, we
need to show that their composition $\dgrm X\circ \dgrm Y$ is a small
functor again. Since $\dgrm X$ is a
weighted colimit of representable functors and small functors are
closed under weighted colimits, 
it suffices to verify that $R^A\circ \dgrm Y$ is a small functor for
any representable functor $R^A$. Then \[(R^A\circ \dgrm Y) (\,\cdot\,) =
R^A(\dgrm Y(\,\cdot\,)) = \map(A, \dgrm Y(\,\cdot\,)) = \dgrm
Y^A. \] Since the category of small
functors is cotensored, $\dgrm Y^A$ is a small functor. 
\end{proof}

\section{A model structure for homotopy functors}\label{homotopy}
A \emph{homotopy functor} is a functor preserving weak equivalences.
In this section, we consider the category $\cal S^{\cal S}$ of small covariant endofunctors  
of simplicial sets and localize the projective model structure 
such that the fibrant functors are exactly the objectwise fibrant homotopy functors. 

Note that small functors are simplicial, hence preserve simplicial homotopy equivalences. 
All simplicial sets are cofibrant, thus every weak equivalence between fibrant simplicial sets is a 
simplicial homotopy equivalence \cite[\S 2, Prop.~5]{Quillen}.
In particular, small functors send weak equivalences of fibrant simplicial sets to weak equivalences. 

We will construct the required localization by the method of Bousfield and Friedlander \cite{BF:gamma}. 
It relies on the existence of a coaugmented functor 
$F\colon \cal S^\cal S \to \cal S^\cal S$, with coaugmentation $\epsilon\colon \Id \to F$. 

\begin{definition}\label{endofunctor}
Given an endofunctor
$F\colon \cat M\to \cat M$ in a model category \cat M equipped with a
coaugmentation $\epsilon\co\Id\to F$ we call a map
$\dgrm{X}\to\dgrm{Y}$ in \cat M an $F$-\emph{equivalence} if it induces a weak equivalence 
$F\dgrm{X}\to F\dgrm{Y}$. A map $\dgrm{X}\to\dgrm{Y}$ is called an 
$F$-\emph{fibration} if it has the right lifting property with respect to all 
projective cofibrations which are also $F$-equivalences. 
\end{definition}

\begin{theorem}[Bousfield-Friedlander]\label{thm:bousfield-machine}
  Suppose $\epsilon \colon \Id\to F$ is a coaugmented endofunctor of a proper model
  category $\cat M$ satisfying the following axioms:
  \begin{description}
  \item[(A.4)] The functor $F$ preserves weak equivalences.
  \item[(A.5)] The maps $\epsilon_{F(A)}, F\epsilon_A \colon F(A) \rightrightarrows FF(A)$ are weak equivalences
    for any object $A\in \cat M$.
  \item[(A.6)] Consider the pullback diagram 
    \diagr{ {W} \ar[r] \ar[d] & {Y} \ar[d]_p \\
        {X} \ar[r]^{f} & {Z}}
      where $p$ is an $F$-fibration and $f$ is an $F$-equivalence.
      Then ${W}\to Y$ is an $F$-equivalence.
  \end{description}
  Then the classes of cofibrations, $F$-equivalences and $F$-fibrations form a
  proper model structure, which is simplicial if $\cat M$ is simplicial.
\end{theorem}

\begin{proof}
  See \cite[A.7]{BF:gamma} and \cite[9.3]{Bou:telescopic}.
\end{proof}
We would like to point out that Theorem~\ref{thm:bousfield-machine} is completely dualizable, and this dual version will be used in Section~\ref{section:homogeneous} to obtain the $n$-homogeneous model structure as a colocalization.

Let ${\rm fib}\colon \cal S\to \cal S$ be a small fibrant replacement functor. For example, 
take ${\rm fib} = \hat{R^\ast} = \hat \Id$, a fibrant replacement of the 
identity functor in the projective model structure on the category of small functors. 
The functor ${\rm fib}$ is equipped with a coaugmentation $\epsilon\colon \Id\to {\rm fib}$.
A concrete example is given by the composition $\mathrm{Sing}\circ |-|$, which is
clearly a simplicial functor. Note that this
functor is small, since it commutes with filtered colimits.
Define $F\colon \cal S^\cal S \to \cal S^\cal S$ by
\begin{equation}\label{F}
      F\dgrm{X} = \dgrm{X}\circ {\rm fib}
\end{equation}
for all $\dgrm{X}\in \cal S^\cal S$. The coaugmentation is given by $\eta_{\dgrm{X}}=\dgrm{X}\circ\epsilon$. 
By \ref{composition is small}, the functor $F\dgrm{X}$ is again small.
Note that a map $\dgrm{X}\to\dgrm{Y}$ is an $F$-equivalence precisely if
$\dgrm{X}(A)\to\dgrm{Y}(A)$ is a weak equivalence for every fibrant
simplicial set $A$. The reason is that $\epsilon_A$ is a simplicial
homotopy equivalence if $A$ is fibrant, and one may conclude with the diagram
\[ \xymatrix @C=5em { 
       \dgrm X (A) \ar[r]^-{\eta_{\dgrm X} (A)=X(\epsilon_A) } \ar[d]_-{f(A)} & (F\dgrm X )(A) =\dgrm{X}\bigl(\mathrm{fib}(A)\bigr) \ar[d]^-{F(f)(A)} \\ 
         \dgrm Y (A) \ar[r]^-{\eta_{\dgrm Y} (A)=Y(\epsilon_A)}  & (F\dgrm X )(A) =\dgrm Y \bigl(\mathrm{fib}(A)\bigr). }\]

\begin{proposition}\label{conditionsA4-6}
The coaugmented functor $\epsilon\colon \Id \to F$ satisfies the axioms 
{\bf (A.4)}, {\bf (A.5)} and {\bf (A.6)}.
\end{proposition}

\begin{proof}
Axiom {\bf(A.4)} holds, because fibrant replacement maps weak equivalences to simplicial homotopy equivalences, and any small functor
preserves these. The maps 
$\epsilon_{{\rm fib}K}, {\rm fib}\,\epsilon_K\co({\rm fib}\circ{\rm fib})(K)\to{\rm fib}(K)$ are weak equivalences 
of fibrant simplicial sets for $K\in\cal{S}$.
Applying the small functor $\dgrm{X}$ preserves these weak equivalences. Hence {\bf (A.5)} holds.
To verify {\bf (A.6)}, observe that pullbacks in the category of small functors are 
computed objectwise. Applying $F$ to the pullback diagram in question thus gives another pullback diagram.
Since $F$-fibrations are in particular objectwise fibrations, $F(p)$ is an objectwise fibration, and $F(f)$ is
an objectwise equivalence. The result follows, because $\cal S$ is right proper. 
\end{proof}

The verification of {\bf (A.6)} shows that $F$-equivalences are closed under base change along
objectwise fibrations.

\begin{theorem}\label{homotopy model structure}
The category of small functors $\cal S^\cal S$ may be equipped with a proper simplicial 
model structure such that weak equivalences are $F$-equivalences, fibrations are $F$-fibrations, and cofibrations are projective cofibrations for $F$ as in {\rm (\ref{F})}.
\end{theorem}

\begin{proof}
  Follows from Theorem~\ref{thm:bousfield-machine} and Proposition~\ref{conditionsA4-6}.
\end{proof}

\begin{definition}
The model structure on $\cal{S}^\cal{S}$ from Theorem \ref{homotopy model structure} will be called the \emph{homotopy model structure}.
\end{definition}

\begin{corollary}\label{cor:char-fib}
A map $\dgrm{X}\to\dgrm{Y}$ is an $F$-fibration if and only if it is an objectwise fibration such that the following square
\diagr{ \dgrm{X} \ar[r]^{\eta_{\dgrm{X}}}\ar[d] & F\dgrm{X} \ar[d] \\
        \dgrm{Y} \ar[r]^{\eta_{\dgrm{Y}}} & F\dgrm{Y}}
is a homotopy pullback square in the projective structure.
In particular, a small functor is $F$-fibrant if and only if it is objectwise fibrant and
a homotopy functor.
\end{corollary}

\begin{proof}
This follows from the characterization given in \cite[Theorem A.7]{BF:gamma}.
\end{proof}

\begin{remark}\label{BF-local vs. Hirschhorn}
\begin{punkt}
     \item
The model structure from \ref{homotopy model structure} coincides with the one constructed in \cite{Lydakis}, provided
one restricts attention to enriched functors from finite pointed simplicial sets to pointed simplicial
sets. What is needed as an essential ingredient is a small enriched fibrant replacement functor -- 
see~\cite[Section 3.3]{DRO:enriched}. To obtain it in our situation, we have to rely on Theorem~\ref{projective structure}
and Lemma~\ref{composition is small}.
    \item
Localization theory implies that a map between $F$-fibrant objects is an 
$F$-equivalence if and only if it is an objectwise equivalence. 
This observation can be seen here directly: obviously a map between homotopy functors is a 
weak equivalence in the homotopy model structure if and only if it is an objectwise equivalence. 
    \item
The coaugmentation $\eta_{\dgrm{X}}\co\dgrm{X}\to F\dgrm{X}$ is not a 
localization in the sense of \cite[3.2.16]{Hirschhorn}, because 
$F\dgrm{X}$ is not required to be fibrant in the projective model structure, 
nor is $\eta_{\dgrm{X}}$ a trivial cofibration in general. 
However, $({\rm fib}\circ F)\dgrm{X}$ is indeed fibrant in the homotopy model structure. 
If we factor the map $\dgrm{X}\to F\dgrm{X}\to({\rm fib}\circ F)\dgrm{X}$ as a 
projective cofibration $\dgrm{X}\to h\dgrm{X}$, followed by an objectwise trivial fibration $h\dgrm{X}\to({\rm fib}\circ F)\dgrm{X}$, 
then we end up with a cofibrant localization in the sense of \cite[3.2.16]{Hirschhorn}. 
The map $\dgrm{X}\to h\dgrm{X}$ has exactly the universal properties that one expects of a 
localization. We point out that if the factorization is not functorial (as it is in our case), then 
$h$ is not functorial in $\dgrm{X}$. But to be clear: $h\dgrm{X}$ is a functor. Here we obtain the following corollary.
\end{punkt}
\end{remark} 

\begin{corollary}
Every small functor may be approximated by a homotopy functor in a universal way up to homotopy. In other words: 
for every small functor $\dgrm X\in \cal S^\cal S$, there exists a functor $h\dgrm X$ and a natural transformation 
$\iota\colon \dgrm X\rarrow h\dgrm X$ such that for every objectwise
fibrant homotopy functor $\dgrm{Y}$ and  every natural transformation $\zeta\colon \dgrm X\to \dgrm Y$, 
there exists a natural transformation $\xi\colon h\dgrm X\to \dgrm Y$, unique up to homotopy, such that $\zeta=\xi\circ\iota$.
\end{corollary}

\begin{proof}
The functor
$h\dgrm X$ is obtained by factorization of the map $\dgrm X\to \ast$ into a 
trivial cofibration followed by a fibration in the homotopy model structure.
The result thus follows from standard model category theory.
\end{proof}

\section{The $n$-excisive structure}
\label{n-excisive structure}

In this section, we localize the homotopy model structure on the category of small endofunctors of \cal S 
such that the fibrant replacement yields the $n$-excisive part of a functor. 
We begin with recalling the relevant definitions from \cite{Goo:calc3}. 

\begin{definition}\label{calclus definitions}
Let $\cal P(\ul{n})$ be the power set of the set $\ul{n}=\{1,...\,,n\}$, equipped with its canonical partial ordering. 
For later use, let $\cal P_0(\ul{n})$ be the complement of $\emptyset$ in $\cal P(\ul n)$. 
An $n$-cubical diagram in $\cal{S}$ is a functor $\cal P(\ul{n})\to\cal{S}$.
A homotopy functor $\dgrm{X}$ is
\begin{itemize}
\item \emph{excisive} if it takes homotopy pushout squares to homotopy pullback squares,
\item \emph{reduced} if $\dgrm{X}(\ast)\simeq \ast$,
\item \emph{linear} if it is both excisive and reduced.
\end{itemize}
A cubical diagram is
\begin{itemize}
\item \emph{strongly homotopy cocartesian} if all of its two-dimensional faces are homotopy pushout squares,
\item \emph{homotopy cartesian} if it is a homotopy limit diagram.
\end{itemize}

A functor $\dgrm{X}$ is said to be $n$\emph{-excisive} if it 
takes strongly homotopy cocartesian $(n+1)$-cubical diagrams to homotopy cartesian diagrams, see \cite[3.1]{Goo:calc3}.
\end{definition}

For an arbitrary homotopy functor $\dgrm{X}$, Goodwillie
\cite[p.~657]{Goo:calc3} constructs an 
$n$-excisive approximation $p_{n,\dgrm{X}}\co\dgrm{X}\to P_n\dgrm{X}$, 
which is natural in $\dgrm{X}$ and universal among all $n$-excisive functors 
under $\dgrm{X}$ \cite[Theorem 1.8]{Goo:calc3}. The homotopy functor $P_n\dgrm{X}$  
is the $n$-excisive part of the Taylor tower of 
$\dgrm{X}$. It is defined as follows: If $U$ is a finite set and $K$ is a simplicial set,
let $K\star U := \hocolim \bigl(K \leftarrow K \times U \rightarrow U\bigr)$.
For any $\dgrm{X}\in \mathcal{S}^\mathcal{S}$, one gets a natural map
\[ t_n\dgrm{X}\colon \dgrm{X}(K) \rightarrow \holim_{\emptyset \neq U \subseteq \{1,\dotsc,n+1\}} \dgrm{X}(K\star U) =: T_n \dgrm{X}(K). \]
Since the homotopy limit of any diagram of simplicial sets is fibrant, $T_n \dgrm{X}$ is objectwise fibrant.
Define $T^\infty_n \dgrm{X}$ to be the sequential colimit of the sequence
\begin{equation}\label{eq:1}\xymatrix@C=4em{ \dgrm{X} \ar[r]^{t_n \dgrm{X}}&  T_n\dgrm{X} \ar[r]^{t_n T_n\dgrm{X}}&  T_n^2 \dgrm{X} \ar[r] &\dotsc} 
\end{equation}
Filtered colimits of simplicial sets preserve weak equivalences, thus $T^\infty_n\dgrm{X}$ is weakly equivalent
to the homotopy colimit of sequence~(\ref{eq:1}). 
Hence $T^\infty_n\dgrm{X}$ is just a simplicial version of Goodwillie's $P_n\dgrm{X}$, in case
$\dgrm{X}$ is a homotopy functor. 

\begin{definition} 
Let $P_n\co\cal{S}^\cal{S}\to\cal{S}^\cal{S}$ be the functor given by
    $$ \dgrm{X}\mapsto P_n\dgrm{X}:=T^\infty_n(\dgrm{X}\circ {\rm fib}). $$
It is a coaugmented functor, with coaugmentation 
$\eta_{n,\dgrm{X}}\colon \dgrm{X} \to \dgrm{X}\circ \mathrm{fib} \to T^\infty_n (\dgrm{X}\circ \mathrm{fib})$.
\end{definition}

Since $T_n\dgrm{X}$ is an objectwise fibrant functor, so is the sequential colimit
$T^\infty_n\dgrm{X}$. If $\dgrm{X}$ preserves weak equivalences, so does $T_n\dgrm{X}$ and
hence also $T^\infty_n\dgrm{X}$. Thus $P_n$ takes values in the category of small objectwise fibrant homotopy functors by construction.

\begin{lemma} \label{properties of P_n}
The functor $P_n$ commutes with finite homotopy limits and 
filtered homotopy colimits. 
\end{lemma}

\begin{proof}
Note that $F$ as in (\ref{F}) commutes with all limits and colimits, because limits and colimits are computed objectwise in the category of small functors, and $F$ is obtained by precomposing with a fixed functor. Fibrations and weak equivalences in the projective model structure are defined objectwise, thus $F$ preserves all homotopy limits. Further, since we work in the category of simplicial sets, filtered colimits are automatically homotopy colimits, thus $F$ preserves these as well. The statement then follows, since $T^\infty_n$ has the claimed properties on the category of homotopy functors by \cite[Prop.~1.7]{Goo:calc3}.
\end{proof}

\begin{proposition}\label{A4-6 again}
The functor $P_n$ satisfies the properties {\bf(A.4)}, {\bf (A.5)} and {\bf (A.6)}.
\end{proposition}

\begin{proof}
Axiom {\bf(A.4)} is fulfilled, because $T_n$ and thus also $T_n^\infty$ preserve objectwise equivalences.
In particular, $P_n$ maps $F$-equivalences to objectwise equivalences.
Axiom {\bf (A.5)} is shown in \cite[proof of 1.8]{Goo:calc3}. 
Axiom {\bf (A.6)} follows directly from the 
fact~\ref{properties of P_n} that 
$P_n$ preserves homotopy pullbacks.
\end{proof}

\begin{definition}
A map $\dgrm{X}\to\dgrm{Y}$ in $\cal{S}^\cal{S}$ is 
\begin{enumerate}
    \item
an $n$-\emph{excisive equivalence} if $P_n\dgrm{X}\to P_n\dgrm{Y}$ is an equivalence in the homotopy model structure.
    \item
an $n$-\emph{excisive fibration} if it has the right lifting property with respect to all projective cofibrations 
which are also $n$-excisive equivalences.
\end{enumerate}
These classes of maps will be called the $n$-\emph{excisive model structure} on $\cal{S}^\cal{S}$.
\end{definition}

The next theorem then follows again from Theorem~\ref{thm:bousfield-machine}
\begin{theorem}
The $n$-excisive structure on $\cal{S}^\cal{S}$ forms a proper simplicial model structure. 
A map $\dgrm{X}\to\dgrm{Y}$ is an $n$-excisive fibration if and only if it is a fibration in the homotopy structure, 
such that the diagram
\diagr{ \dgrm{X} \ar[r]^-{\eta_{n,\dgrm{X}}}\ar[d] & P_n\dgrm{X} \ar[d] \\
        \dgrm{Y} \ar[r]^-{\eta_{n,\dgrm{Y}}} & P_n\dgrm{Y} }
is a homotopy pullback square in the homotopy structure. Fibrant objects are exactly the 
objectwise fibrant $n$-excisive homotopy functors.
\end{theorem}

\begin{remark}\label{remark:strong-right-proper}
Since functors of the form $P_n\dgrm{X}$ are homotopy functors by construction, $\dgrm{X}\to\dgrm{Y}$ is an $n$-excisive equivalence 
if and only if $P_n\dgrm{X}\to P_n\dgrm{Y}$ is an objectwise equivalence using remark \ref{BF-local vs. Hirschhorn}(ii).
Note also that $n$-excisive equivalences are closed under base change
along objectwise fibrations. This follows, since $P_n$ preserves homotopy pullbacks
in the projective model structure.
\end{remark}

\section{Homotopy theory of spectrum-valued functors}
\label{section:spectrum-valued}

In this section, we introduce a model category that describes homotopy
theory of small functors with values in spectra. First of all we have to
give a definition of small spectrum-valued functors. To streamline the exposition, we will use 
the category of pointed spaces $\cal{S}_\ast$ as our underlying symmetric monoidal category. 
Note that the results obtained in Sections~\ref{prelim},~\ref{homotopy} and~\ref{n-excisive structure} 
go through for the category $\cal{S}_*^{\cal{S}_*}$ 
of small endofunctors of pointed spaces after making the canonical changes. 
Let $\Spec$ denote the category of spectra in the sense of
Bousfield-Friedlander \cite{BF:gamma}. The category of spectra in a pointed simplicial model category $\cal{M}$ 
is denoted $\Spec(\cal{M})$.

\begin{remark}\label{alternative}
Note that every small functor $\dgrm X\in \cal{S}_*^{\cal{S}_*}$ preserves the initial object $\dgrm X(\ast)=\ast$, since it is a weighted colimit of representable functors. In other words all small functors in $\cal{S}_*^{\cal{S}_*}$ are reduced.

An alternative way to consider the homotopy theory of functors from
$\cal S_\ast$ to $\cal S_\ast$ is to look at the category of small
\emph{simplicial} functors from $\cal S_\ast$ to $\cal S$ and then to
form the category under the constant functor $\ast$. The main
difference between these approaches is, that the resulting category of
functors will be enriched over simplicial sets, rather then pointed
simplicial sets. As a consequence, elements of $\ast\overcat \cal S^{\cal S_\ast}$ include also non-reduced functors.
\end{remark}

\begin{definition}
An object in the category $\Spec(\cal{M})$ is a sequence $(X_0,X_1,...)$ of objects in $\cal{M}$, together with bonding maps
    $$ \Sigma X_n \to X_{n+1} $$
for $n\ge 0$, where $\Sigma X_n:= X_n\otimes \Delta^1/\partial\Delta^1$.
\end{definition}

\begin{definition}
A functor from $\cal{S}_*$ to $\Spec$ is \emph{small} if it is the enriched left Kan extension of a functor defined on a small subcategory of $\cal{S}_*$. 
\end{definition}

For each $n\ge 0$ let ${\rm Ev}_n\co\Spec\to\cal{S}_*$ denote the functor taking a spectrum $X=(X_0,X_1,...)$ to its $n$-th level $X_n$.

\begin{lemma} \label{small=levelwise small}
A functor $\dgrm{X}\co\cal{S}_*\to\Spec$ is small if and only if it is 
levelwise small, i.e. if ${\rm Ev}_n\circ \dgrm{X}\co\cal{S}_*\to\cal{S}_*$ is small for each $n\ge 0$.
\end{lemma}

\begin{proof}
The evaluation functors ${\rm Ev}_n$ are simplicial and have enriched right adjoints, which therefore commute with enriched left Kan extensions.
\end{proof}

\begin{lemma}\label{evident}
The evident functors give an isomorphism $\Spec(\cal{S}_*^{\cal{S}_*})\cong\Spec^{\cal{S}_*}$ of categories.
\end{lemma}

\begin{proof}
This follows directly from Lemma \ref{small=levelwise small}.
\end{proof}

Using Lemma \ref{evident} we will identify the categories $\Spec(\cal{S}_*^{\cal{S}_*})$ and $\Spec^{\cal{S}_*}$.
This shows in particular that the category $\Spec^{\cal{S}_*}$ is complete. 
Now we want to lift the projective model structure, where a weak equivalence 
is given objectwise, to the spectrum-valued setting.
Our strategy is the following: We take the projective model structure on 
$\cal{S}_*^{\cal{S}_*}$ and then consider spectrum objects over this category. 
Using results from \cite{Schwede:cotangent} we obtain a model structure on 
$\Spec(\cal{S}_*^{\cal{S}_*})$, which is the desired one.


In \cite{Schwede:cotangent} the stable model structure on spectra is obtained analogously as in \cite{BF:gamma}, 
but the construction $Q$ used there is adaptable to more general situations. Lemma 1.3.2. of 
\cite{Schwede:cotangent} lists the properties which have to be satisfied in order to make the machinery work. 
Although in our case the underlying model structure on $\cal{S}_*^{\cal{S}_*}$ is not cofibrantly generated, 
we are still able to prove the required statements. The reason is that our model category is class-cofibrantly generated. 
Here is the adapted version of the relevant part (a) of Lemma 1.3.2. of \cite{Schwede:cotangent}.

\begin{lemma}\label{lem:finitely-generated}
Let $X\to Y$ be a termwise $($trivial$)$ fibration between sequences 
in the category $\cal{S}_*^{\cal{S}_*}$. Then the induced map $\colim X\to\colim Y$ is a $($trivial$)$ fibration. In particular, 
sequential colimits preserve weak equivalences.
\end{lemma}

\begin{proof}
The proof for the case of fibration and trivial fibration is literally the same except that 
one uses the different test classes $I$ or $J$ from (\ref{IandJ}). Since source and target 
of the generating classes $I$ and $J$ are small, we get the following liftings
\diagr{ R^A\otimes K \ar[r]\ar[d]_{i} & X_k \ar[r]\ar@{->>}[d]^{(\simeq)} & \colim X \ar[d] \\
        R^A\otimes L \ar[r]\ar@{.>}[ur] & Y_k \ar[r] & \colim Y}
where $i$ is either in $I$ or $J$. This proves the statement.
\end{proof}

For the definition of the coaugmented functor $Q\co\Spec(\cal{S}_*^{\cal{S}_*})\to\Spec(\cal{S}_*^{\cal{S}_*})$ we refer to \cite[p. 93]{Schwede:cotangent}. 
Note, however, that one may use a simpler construction instead, which avoids factoring the bonding maps -- see \cite{Hovey:spectra}.
Given a spectrum $(X_0,X_1,\dotsc)$ with bonding maps adjoint to $\sigma_n\colon X_n\to X_{n+1}^{S^1}$, 
let $\Omega_\mathrm{fake}X[1]$ denote the spectrum with
terms $\bigl(\Omega_\mathrm{fake}X[1]\bigr)_n = X_{n+1}^{S^1}$ and bonding maps adjoint to 
\[\sigma_{n+1}^{S^1}\colon \bigl(\Omega_\mathrm{fake}X[1]\bigr)_n = X_{n+1}^{S^1} \to X_{n+2}^{S^1} = \Omega_\mathrm{fake}X[1]_{n+1}.\] 
The adjoints of the bonding maps define a natural map $X\to \Omega_\mathrm{fake}X[1]$ of spectra.
Given a spectrum $\dgrm{X}$ of small functors, define $Q'(\dgrm{X})$ as the colimit of the sequence
\[ \xymatrix{ \mathrm{fib}\circ \dgrm{X} \ar[r] & \Omega_\mathrm{fake}\bigl(\mathrm{fib}\circ \dgrm{X}\bigr)[1]  \ar[r] & \dotsc}.\]
Then Lemma~\ref{lem:finitely-generated} shows that the canonical map $Q\dgrm{X} \to Q'\dgrm{X}$ is an objectwise
weak equivalence in each level. 
For each $K$ in $\cal{S}_*$ the spectrum $(Q\dgrm{X})(K)$ is weakly equivalent to the usual 
$\Omega$-spectrum $Q(\dgrm{X}(K))$ in the Bousfield-Friedlander sense.
 
\begin{definition}
A map $\dgrm{X}\to\dgrm{Y}$ in $\Spec(\cal{S}_*^{\cal{S}_*})$ will be called 
\begin{punkt}
    \item
a \emph{projective cofibration} if the map $\dgrm{X}_0\to\dgrm{Y}_0$ and for each $n\ge 0$ 
the maps $\dgrm{X}_n\sqcup_{\dgrm{X}_{n-1}}\dgrm{Y}_{n-1}\to\dgrm{Y}_n$ are projective cofibrations.
    \item
a \emph{stable objectwise equivalence} if for all $n\ge 0$ the maps $Q\dgrm{X}_n\to Q\dgrm{Y}_n$ are objectwise equivalences.
    \item
a \emph{stable objectwise fibration} if for all $n\ge 0$ the maps $Q\dgrm{X}_n\to Q\dgrm{Y}_n$ are objectwise fibrations and the squares
\diagr{ \dgrm{X}_n \ar[r]\ar[d] & Q\dgrm{X}_n \ar[d] \\
        \dgrm{Y}_n \ar[r] & Q\dgrm{Y}_n }
are homotopy pullback squares in the projective structure.
\end{punkt}
We call these classes of maps the \emph{stable projective model structure} on $\Spec(\cal{S}_*^{\cal{S}_*})$.
\end{definition}

\begin{proposition}\label{prop:stable-projective}
The stable projective model structure on $\Spec(\cal{S}_*^{\cal{S}_*})\cong\Spec^{\cal{S}_*}$ is a simplicial proper model structure.
\end{proposition}

The proof of Proposition~\ref{prop:stable-projective} is as the proof of~\cite[Prop.~2.1.5]{Schwede:cotangent}. 
We do not claim that the stable projective model structure has
 -- or has not -- functorial factorization. 

Since all functors in $\Spec^{\cal{S}_*}$ are simplicial, they preserve simplicial
homotopies and therefore map weak equivalences between fibrant spaces
to weak equivalences. The same method as in
Section~\ref{homotopy} thus provides the homotopy model structure on
$\Spec^{\cal{S}_*}$.

\begin{definition}
A map $\dgrm{X}\to\dgrm{Y}$ in $\Spec^{\cal{S}_*}$ is called 
\begin{punkt}
    \item
a \emph{stable equivalence in the homotopy structure} if $\dgrm{X}(K)\to\dgrm{Y}(K)$ is a stable equivalence of spectra for all fibrant spaces $K$.
    \item
a \emph{stable fibration in the homotopy structure} if $\dgrm{X}(K)\to\dgrm{Y}(K)$ is a stable fibration for all spaces $K$ and the square
\diagr{ \dgrm{X} \ar[r]\ar[d] & \dgrm{X}\circ {\rm fib}=F\dgrm{X} \ar[d] \\
        \dgrm{Y} \ar[r] & \dgrm{Y}\circ {\rm fib}=F\dgrm{Y} }
is a homotopy pullback square in the stable projective structure. 
\end{punkt}
We call these classes of maps the stable homotopy model structure on $\Spec(\cal{S}_*^{\cal{S}_*})$.
\end{definition}

As in Section~\ref{homotopy} we obtain the following theorem.
\begin{theorem}
The stable homotopy model structure on $\Spec^{\cal{S}_*}$ is a simplicial proper model structure. 
A functor in $\Spec^{\cal{S}_*}$ is a homotopy functor if and only if it is weakly equivalent in the 
stable projective structure to a fibrant object in the stable homotopy structure.
\end{theorem}

There are different characterizations of weak equivalences, here we give one.
\begin{lemma}
A map $\dgrm{X}\to\dgrm{Y}$ is a weak equivalence in the stable homotopy structure if 
and only if for each $n\ge 0$ the maps $Q\dgrm{X}^n\to Q\dgrm{Y}^n$ are weak equivalences in the homotopy structure on $\cal{S}_*^{\cal{S}_*}$.
\end{lemma}

\begin{proof}
This follows from the natural equivalence $Q(\dgrm{X}\circ {\rm fib})\cong (Q\dgrm{X})\circ {\rm fib}$. 
\end{proof}

Since $\holim$ and $\hocolim$ are simplicial functors, so is the functor $T^\infty_n$. Hence $P_n$ has a natural extension
to functors with values in spectra. Localizing along the coaugmented functor $\Id \to P_n$ then
yields the stable $n$-excisive model structure on $\Spec^{\cal{S}_*}$.
The following property of $P_n$, which will be relevant later on, may be deduced from~\cite[Prop.~1.7]{Goo:calc3}.

\begin{lemma}\label{lem:pn-spectra}
  The functor $P_n\colon \Spec^{\cal{S}_*} \to \Spec^{\cal{S}_*}$ commutes with homotopy colimits.
\end{lemma}

\begin{theorem}
The stable $n$-excisive model structure on $\Spec^{\cal{S}_*}$ is a
simplicial proper model structure. A functor in $\Spec^{\cal{S}_*}$ is
an $n$-excisive homotopy functor if and only if it is weakly
equivalent in the stable projective structure to a fibrant object in
the stable $n$-excisive structure.
\end{theorem}

\section{The Taylor tower and homogeneous functors}
\label{section:Taylor-tower}

In this section, the homotopy theory under consideration is the $n$-excisive model structure on
the categories $\cal{S}_\ast^{\cal{S}_\ast}$ and $\Spec^{\cal{S}_*}$. The existence of basepoints in the target category is required 
for taking certain homotopy fibers. A pointed source category is not required here, and everything in sections \ref{section:Taylor-tower} and \ref{section:homogeneous} applies to small functors initiating in \cal{S}. For expositional reasons we have chosen $\cal{S}_*$, since in section \ref{section:F-equiv} we will switch to a pointed source category.

The fibrant objects in the $n$-excisive structure are the objectwise fibrant $n$-excisive homotopy functors. 
Since $(n-1)$-excisive functors are also $n$-excisive, there is a map $P_n\dgrm{X}\to P_{n-1}\dgrm{X}$ under $\dgrm{X}$ in the homotopy category. 
By the results from the previous section, this map is unique up to simplicial homotopy.
Goodwillie gives a model for this map 
for homotopy functors of topological spaces \cite[p. 664]{Goo:calc3}.
One may immediately obtain the same natural map 
    $$ q_{n,\dgrm{X}}\co P_n\dgrm{X}\to P_{n-1}\dgrm{X} .$$
in any of the categories under consideration. These maps fit into a tower under 
$\dgrm{X}$, which is called the Taylor tower of $\dgrm{X}$. The fibers of this tower are of special interest. Let us 
recall a definition first. 

\begin{definition}
A functor $\dgrm{X}$ is called $n$-\emph{reduced} if $\dgrm{X}$ is 
weakly contractible in $(n-1)$-excisive structure, i.e. 
$P_{n-1}\dgrm{X}\simeq\ast$ in the homotopy structure. 
A functor is called $n$-\emph{homogeneous} if it is $n$-reduced and $n$-excisive.
\end{definition}

To introduce the homogeneous part $D_n\dgrm{X}$ of a 
small functor $\dgrm{X}$ we consider fibers and homotopy fibers. 
Recall that the simplicial path object of an object $Z$ in a pointed simplicial model category is given by
    $$ WZ := Z^{\Delta^1}\times_{(Z\times Z)}(\ast\times Z),$$
where the map $Z^{\Delta^1}\to Z\times Z$ is induced by 
$d_0\vee d_1\co \Delta^0\vee\Delta^0\to\Delta^1$. 
The projection ${\rm pr}_2\co Z\times Z\to Z$ induces a map $WZ\to Z$. If $Z$ is fibrant, this map is a fibration.
Note that $WZ$ is simplicially contractible.

\begin{definition}\label{defn:Dn}
If $\dgrm{X}$ is a small functor, a new small functor $D_n\dgrm{X}$ is defined by the following pullback square:
\diag{ D_n\dgrm{X} \ar[r]\ar[d]_{d_{n,\dgrm{X}}} & W(P_{n-1}\dgrm{X}) \ar[d] \\
       P_n\dgrm{X} \ar[r]_{q_{n,\dgrm{X}}} & P_{n-1}\dgrm{X}  }{def. of D_n}
We call $D_n\dgrm{X}$ the $n$-\emph{homogeneous part} of $\dgrm{X}$.
\end{definition}

\begin{remark}
The map $q_n\dgrm{X}$ is an equivalence in the $(n-1)$-excisive structure, 
and therefore $D_n\dgrm{X}$ is $(n-1)$-excisively contractible, hence $n$-reduced. 
Since $d_{n,\dgrm{X}}$ is the base change of an $n$-excisive fibration, 
$D_n\dgrm{X}$ is $n$-excisively fibrant, and thus $n$-homogeneous. 
We also point out that the square (\ref{def. of D_n}) is a homotopy pullback square in the 
following model structures: in the projective, the homotopy, and the $(n-1)$-excisive structure. 
\end{remark}

We will need the following properties, which are given in \cite[Prop. 1.18]{Goo:calc3}.
\begin{proposition}\label{prop:dn-commutes}
The functor $D_n\co\cal{S}_*^{\cal{S}_*}\to\cal{S}_*^{\cal{S}_*}$ commutes with 
finite homotopy limits and filtered homotopy colimits in the projective and the homotopy model structure. 
The functor $D_n\co\Spec^{\cal{S}_*}\to\Spec^{\cal{S}_*}$ commutes with finite homotopy limits and 
all homotopy colimits in the projective and the homotopy model structure.
\end{proposition}

\section{The $n$-homogeneous structure}
\label{section:homogeneous}

In this section, we construct the $n$-homogeneous model structure on $\Spec^{\cal{S}_*}$
via a colocalization process, which involves the dual of Theorem~\ref{thm:bousfield-machine}. 
We only claim this structure for the spectrum-valued case as will become apparent in Lemma \ref{left properness for M_n}.
The homotopy types correspond bijectively to $n$-homogeneous spectrum-valued functors. 
This model structure classifies \emph{all} $n$-homogeneous functors as homotopy types 
(cf. \cite[Remark 4.13]{Kuhn:overview}). We will give an interpretation of Goodwillie's 
classification of \emph{finitary} homogeneous functors in terms of a Quillen equivalence 
between model categories in Section~\ref{section:F-equiv}.

  
The situation is similar to the case of the ordinary Postnikov tower of spaces. 
One obtains each Postnikov stage as a fibrant replacement by localizing with respect 
to $S^n$, which kills all homotopy above degree $n$. One can also colocalize with respect 
to $S^{n-1}$. This time the cofibrant replacement is given by the connected covers; the 
homotopy below degree $n-1$ is killed.
Here we are going to colocalize with respect to the $n$-reduced part of a functor.

\begin{definition}
For each small functor $\dgrm{X}$, let $M_n\dgrm{X}$ be defined by the following pullback square:
\diagr{ M_n\dgrm{X} \ar[r]\ar[d]_{m_{n,\dgrm{X}}} & W(P_{n-1}\dgrm{X}) \ar[d] \\
        \dgrm{X} \ar[r]_-{p_{n-1,\dgrm{X}}} & P_{n-1}\dgrm{X}  }
The augmented functor $M_n\co\cal{S}_*^{\cal{S}_*}\to\cal{S}_*^{\cal{S}_*}$ is called the $n$-\emph{reduced part} of $\dgrm{X}$. 
\end{definition}

\begin{remark}
The object $P_{n-1}\dgrm{X}$ is fibrant in the $(n-1)$-excisive model structure, hence
$WP_{n-1}\dgrm{X}\to P_{n-1}\dgrm{X}$ is an $(n-1)$-excisive fibration.
By right properness, it follows that $M_n\dgrm{X}$ is the homotopy pullback of 
$\dgrm{X}\to P_{n-1}\dgrm{X}\leftarrow W(P_{n-1}\dgrm{X})$ in the $(n-1)$-excisive model structure,
thus also in the homotopy and the projective model structure. 

Since the map $p_{n-1}\dgrm{X}$ is an $(n-1)$-excisive equivalence, 
the functor $M_n\dgrm{X}$ is weakly contractible in the $(n-1)$-excisive structure, 
and therefore $n$-reduced. For each $\dgrm{X}$ we have a square
\diagr{ M_n\dgrm{X} \ar[r]\ar[d]_{m_{n,\dgrm{X}}} & D_n\dgrm{X} \ar[d]^{d_{n,\dgrm{X}}} \\
        \dgrm{X} \ar[r]_{p_n\dgrm{X}} & P_n\dgrm{X} } 
which is a pullback as well as a homotopy pullback square in the projective, homotopy or $(n-1)$-excisive structure. 
The construction  $M_n$ preserves homotopy pullbacks, since it is the homotopy fiber of functors preserving homotopy pullbacks.
Of course, $M_n\dgrm{X}$ is not a homotopy functor unless $\dgrm{X}$ is one.
\end{remark}

To colocalize along the functor $M_n$, we have to prove that the axioms dual to 
the Bousfield-Friedlander axioms in Theorem~\ref{thm:bousfield-machine} are satisfied. 
As a starting model structure for the colocalization, one may use either the $n$-excisive structure, or the homotopy structure. 
In the first case we obtain the $n$-homogeneous structure \ref{n-homogeneous structure}; in the second case the resulting model structure is 
called the $n$-reduced structure \ref{n-reduced structure}. The statements \ref{A4-dual}, \ref{A5-dual} and 
\ref{left properness for M_n} cover both cases. The proof of the left properness condition 
{\bf (A.6)}${}^\op$ uses that $D_n$ commutes with homotopy pushouts.
This explains the restriction to spectrum-valued functors.
 
\begin{lemma} \label{A4-dual}
The functor $M_n$ satisfies {\bf (A.4)}${}^\op$.
\end{lemma}

\begin{proof}
Since $M_n\dgrm{X}$ is defined as the homotopy pullback in the homotopy model structure of the 
diagram $\dgrm{X} \to P_{n-1}\dgrm{X} \leftarrow \ast$, any weak equivalence $\dgrm{X}\to \dgrm{Y}$
in the homotopy structure induces a weak equivalence $M_n\dgrm{X} \to M_n\dgrm{Y}$ in the homotopy structure. 
Suppose now that $f\colon \dgrm{X}\to \dgrm{Y}$ is an $n$-excisive equivalence. Since $n$-excisive equivalences
are closed under base change along objectwise fibrations by Remark~\ref{remark:strong-right-proper}, 
the map $M_n\dgrm{X}\to D_n\dgrm{X}$ is an $n$-excisive equivalence. It remains
to check that $D_n$ preserves $n$-excisive equivalences. Since $n$-excisive equivalences are in particular
$(n-1)$-excisive equivalences, $P_n(f)$, $P_{n-1}(f)$ and $WP_{n-1}(f)$ are objectwise equivalences.
The result follows from Proposition~\ref{prop:dn-commutes}, 
because $WP_{n-1}\dgrm{X} \to P_{n-1}\dgrm{X}$ is in particular an objectwise fibration.
\end{proof}

\begin{lemma} \label{A5-dual}
The maps $m_{n,M_n\dgrm{X}}$ and
$M_nm_{n,\dgrm{X}}\co M_nM_n\dgrm{X}\to M_n\dgrm{X}$ are objectwise equivalences.
Hence $M_n \to \Id$ satisfies {\bf (A.5)}${}^\op$.
\end{lemma}

\begin{proof}
The map $m_{n,M_n\dgrm{X}}$ is an objectwise equivalence, because it is the base change
of the objectwise acyclic fibration $WP_{n-1}M_n\dgrm{X}\to P_{n-1}M_n\dgrm{X}$. 
Since $P_{n-1}\dgrm{X}\to P_{n-1}P_{n-1}\dgrm{X}$ is an objectwise equivalence,
so is the base change $M_nP_{n-1}\dgrm{X}\to WP_{n-1}P_{n-1}\dgrm{X}$. Thus
$M_n WP_{n-1}\dgrm{X}\to M_nP_{n-1}\dgrm{X}$ is an objectwise equivalence.
The result follows, since $M_n$ preserves homotopy pullbacks.  
\end{proof}

\begin{lemma} \label{left properness for M_n}
The augmentation $M_n\to \Id$ satisfies {\bf (A.6)}${}^\op$.
\end{lemma}

\begin{proof}
Consider the following diagram
   $$\xy
   \xymatrix"*"@=17pt{ \dgrm{A} \ar[dd]_j \ar[rr] & & \dgrm{X} \ar@{->>}[dd] \\
                                 & &                \\
                    \dgrm{B} \ar[rr] & &  \dgrm{Y}  }
   \POS(9,8)
   \xymatrix@=11pt{ M_n\dgrm{A} \ar[rr]\ar'[d] [dd]\ar["*"]  & & M_n\dgrm{X} \ar[dd]\ar["*"]  \\
                                 & &                \\
              M_n\dgrm{B} \ar'[r] [rr]\ar["*"]   & & M_n\dgrm{Y} \ar["*"]   }
\endxy   $$
where $j$ is a cofibration 
and the front square is a pushout square. The homotopy model structure is left proper, thus the
front square is a homotopy pushout square.
Since the spectrum-valued functor $M_n$ preserves homotopy pushouts in the homotopy model 
structure, the back square is a homotopy pushout in the homotopy model structure. 
This already proves {\bf (A.6)}${}^\op$ in the case of the homotopy model structure: 
If the map $M_nA\to M_nX$ is an $F$-equivalence for $F$ as in (\ref{F}), then so is $M_nB\to M_nY$. 
 
To prove {\bf (A.6)}${}^\op$ in the case of the $n$-excisive model structure, 
let $M_n\dgrm{A}\to M_n\dgrm{X}$ be an $n$-excisive equivalence. We have to show that the map $M_n\dgrm{B}\to M_n\dgrm{Y}$ is 
also an $n$-excisive equivalence. In general, a map $\dgrm{S}\to\dgrm{T}$ is an $n$-excisive 
equivalence if and only if $P_n\dgrm{S}\to P_n\dgrm{T}$ is an objectwise equivalence. 
So we apply $P_n$ to the back square. Using the objectwise equivalence $P_nM_n\dgrm{X} \to D_n\dgrm{X}$, 
we obtain the square
\diag{ D_n\dgrm{A} \ar[r]^\simeq\ar[d] & D_n\dgrm{X} \ar[d] \\
       D_n\dgrm{B} \ar[r] & D_n\dgrm{Y}  }{D_n}
where $D_n\dgrm{A}\to D_n\dgrm{X}$ is an objectwise equivalence. 
Since the spectrum-valued functor $D_n$ preserves homotopy pushouts by~\ref{prop:dn-commutes}, 
the square (\ref{D_n}) is a homotopy pushout square in the projective model structure. 
In particular, $D_n\dgrm{B}\to D_n\dgrm{Y}$ is an objectwise equivalence, 
thus $M_n\dgrm{B}\to M_n\dgrm{Y}$ is an $n$-excisive equivalence.
\end{proof}

\begin{definition} 
A map $f\co\dgrm{X}\to\dgrm{Y}$ in $\Spec^{\cal{S}_*}$
is an $n$-\emph{homogeneous equivalence} if
    $$ D_n{\dgrm{X}}\to D_n{\dgrm{Y}} $$
is an $n$-excisive equivalence. The map $f$ is an $n$-\emph{homogeneous cofibration} if it
has the left lifting property with respect to all $n$-excisive fibrations which are
$n$-homogeneous equivalences. These classes form the \emph{$n$-homogeneous structure} on~$\Spec^{\cal S_\ast}$.
\end{definition}

Using remark \ref{BF-local vs. Hirschhorn}(ii) repeatedly, 
$\dgrm{X}\to\dgrm{Y}$ is an $n$-homogeneous equivalence if 
and only if $D_n{\dgrm{X}}\to D_n{\dgrm{Y}}$ is an objectwise equivalence.
Observe also that $M_n(f)$ is an $n$-excisive equivalence if and only if $D_n(f)$ is
an $n$-excisive equivalence.
Since the axioms $\mathbf{(A.4)^\op}$, $\mathbf{(A.5)^\op}$ and $\mathbf{(A.6)^\op}$ hold,
we obtain the following statement from the dual of Theorem~\ref{thm:bousfield-machine}.
\begin{theorem} \label{n-homogeneous structure}
On the category $\Spec^{\cal S_\ast}$ the $n$-homogeneous structure is a proper simplicial
model structure. The fibrant and cofibrant objects are exactly the projectively
cofibrant $n$-homogeneous homotopy functors having values in stably fibrant spectra. 
In particular, the homotopy types correspond 
bijectively to the homotopy types of $n$-homogeneous functors from $\cal{S}_*$ to $\Spec$.
\end{theorem}


For a small functor $\dgrm{X}$ the object $P_n\dgrm{X}$ is not exactly the localization of 
$\dgrm{X}$ in the sense of \cite[3.2.16]{Hirschhorn}, as explained in remark \ref{BF-local vs. Hirschhorn}(iii). 
But $P_n\dgrm{X}$ is not far away from that; it is weakly equivalent to the localization in the underlying model 
structure, here the homotopy structure. The same is true for $D_n\dgrm{X}$: The maps $\dgrm{X}\to P_n\dgrm{X}\leftarrow D_n\dgrm{X}$ 
are not a fibrant approximation followed by a cofibrant approximation, but $D_n\dgrm{X}$ is weakly equivalent in the homotopy 
structure to a fibrant and cofibrant replacement of $\dgrm{X}$ in the $n$-homogeneous structure. In fact, since 
both functors $D_n\dgrm{X}$ and the replacement of $\dgrm{X}$ are homotopy functors, they are even weakly 
equivalent in the projective structure on $\cal{S}_*^{\cal{S}_*}$.

Finally it is worth remarking that we can colocalize along the functor $M_n$ 
starting directly from the homotopy structure.

\begin{theorem} \label{n-reduced structure}
The category $\Spec^{\cal{S}_*}$ may be equipped with the $n$-reduced
model structure. The resulting model category is simplicial and
proper. The cofibrant objects are exactly the projectively cofibrant
$n$-reduced functors. 
\end{theorem}

\section{Spectra with $\Sigma_n$-action and $n$-homogeneous functors}
\label{section:F-equiv}

The goal of this section is to connect the homotopy theory of small
spectrum-valued $n$-homogeneous functors with the homotopy theory of
spectra with $\Sigma_n$-action. We interpret the $n$-th derivative at
$\ast$ as a part of a Quillen pair between these categories. In
Section~\ref{section:finitary} we will show that this pair induces a
Quillen equivalence when applied to the category of small functors
with the finitary homogeneous model structure. Altogether,
this may be viewed as a strengthening of Goodwillie's result
\cite{Goo:calc3}. 

\begin{remark} \label{naive}
\begin{punkt}
   \item
We consider spectra with $\Sigma_n$-action, that is, presheaves on the group $\Sigma_n$ (considered as a category with one object) 
with values in spectra. We equip spectra with the Bousfield-Fried\-lan\-der model structure and take the projective 
model structure on presheaves over it. Thus, weak equivalences or fibrations are given by weak equivalences or 
fibrations of the underlying spectra. It is sometimes called the naive
equivariant model structure. We will denote it by $\Spec^{\Sigma_n}$.
   \item
More generally, one can endow $\Sigma_n$-objects in any cofibrantly
generated model category with the projective model structure
\cite{Hirschhorn}. It is easy to check that this result holds for any
class-cofibrantly generated model category. Since the category
$\cal{S_*}^{\cal{S}_*}$  with the projective model structure is
class-cofibrantly generated \cite{Chorny-Dwyer}, the category of small
functors with $\Sigma_n$-action may be given the naive equivariant model structure.
   \item
In order to exhibit certain maps and diagrams in Lemmas~\ref{cross
  effect} and \ref{lem:stable-cross-effect} we need a basepoint in the
source category. In the same way Goodwillie \cite[p. 676]{Goo:calc3}
considers cross effects only for pointed source categories. 
\end{punkt}
\end{remark}

Let us start to define the Quillen pair $\lambda_n\co\Spec^{\Sigma_n}\leftrightarrows\Spec^{\cal{S}_*}:\!\rho_n$.
First we need a homotopy invariant version of the objectwise smash product in $\cal{S}_*^{\cal{S}_*}$. 
The objectwise smash product $\dgrm{X}\wedge \dgrm{Y}$ is given by the quotient of
the canonical map $\dgrm{X}\vee \dgrm{Y} \hookrightarrow \dgrm{X}\times \dgrm{Y}$. Since this map
is not a projective cofibration in general, the objectwise smash product might fail to be cofibrant -- 
even for representable functors like $\id\cong R^{S^0}$. 
For $K_1, ...\, ,K_n\aus\cal{S}_*$ let 
    $$ \left(\bigwedge_{i=1}^nR^{K_i}\right)_{\rm cof} \to \bigwedge_{i=1}^nR^{K_i} $$
be a projective cofibrant replacement. 
We can relate the smash product to the $n$-th cross effect as defined in \cite[p.~676]{Goo:calc3} or \cite[5.8]{Kuhn:overview}.
Recall from \ref{calclus definitions} that $\cal P_0(\ul{n})=\cal P(\ul{n})-\emptyset$.

\begin{lemma} \label{cross effect}
\begin{punkt}
   \item
For any small functor $\dgrm{X}\colon \mathcal{S}_\ast \rightarrow \mathcal{S}_\ast$ there is a natural isomorphism:
\begin{align*}
    \map\left(\bigwedge\limits_{i=1}^nR^{K_i},\dgrm{X}\right) \cong {\rm
    fib}\left[\dgrm{X}\left(\bigvee_{i=1}^nK_i\right)\to\lim_{T\in \cal P_0(\ul{n})}\dgrm{X}\left(\bigvee_{\ul{n}-T}K_i\right) \right] 
\end{align*}
   \item
If $\dgrm{X}$ is objectwise fibrant, there is a natural objectwise weak equivalence:
\begin{align*}
    \map\left(\left(\bigwedge\limits_{i=1}^nR^{K_i}\right)_{\rm cof},\dgrm{X}\right)
                   &\simeq {\rm
                   hofib}\left[\dgrm{X}\left(\bigvee_{i=1}^nK_i\right)\to\holim_{T\in
                   \cal P_0(\ul{n})}\dgrm{X}\left(\bigvee_{\ul{n}-T}K_i\right) \right] \\ 
                   & \cong{\rm cr}_n\dgrm{X}(K_1,...\, ,K_n) 
\end{align*}
\end{punkt}
\end{lemma}

\begin{proof}
  There is an enriched Yoneda isomorphism
  \[ \map(R^K,\dgrm{X}) \cong X(K).\]
  Part (i) then follows by adjunction from the representation of an iterated smash product as the pushout diagram
  \diagr{ \displaystyle{\colim_{T\in \cal P_0(\ul{n})}}\prod_{i\in\ul{n}-T}R^{K_i} \ar[r]\ar[d] & \prod_{i=1}^n R^{K_i} \ar[d] \\
        \ast \ar[r] & \bigwedge_{i=1}^n R^{K_i} }
  Part (ii) follows from (i), because the source is cofibrant and the target is fibrant. 
\end{proof}

The analog of Lemma~\ref{cross effect} for small functors $\dgrm{X}\colon \mathcal{S}_\ast \rightarrow \Spec$ 
is obtained as follows: For $\dgrm{K} \in \cal{S}_*^{\cal{S}_*}$ and $\dgrm{X} \in \Spec^{\cal{S}_*}$ let
$\mathrm{spt}(\dgrm{K},\dgrm{X})$ be the spectrum whose $k$-th term is
\[ \mathrm{spt}(\dgrm{K},\dgrm{X})_k \colon = \map(\dgrm{K},\mathrm{Ev}_k\circ \dgrm{X}).\]
One obtains bonding maps for $\mathrm{spt}(\dgrm{K},\dgrm{X})$, because $\cal{S}_*^{\cal{S}_*}$ is enriched, tensored and cotensored
over $\mathcal{S}_\ast$. There is an enriched Yoneda isomorphism 
    $$ \mathrm{spt}(R^K,\dgrm{X})\cong\dgrm{X}(K).$$

\begin{lemma} \label{lem:stable-cross-effect}
\begin{punkt}
   \item
For any small functor $\dgrm{X}\colon \mathcal{S}_\ast \rightarrow \Spec$ there is a natural isomorphism:
\begin{align*}
    \mathrm{spt}\left(\bigwedge\limits_{i=1}^nR^{K_i},\dgrm{X}\right) \cong {\rm
    fib}\left[\dgrm{X}\left(\bigvee_{i=1}^nK_i\right)\to\lim_{T\in \cal P_0(\ul{n})}\dgrm{X}\left(\bigvee_{\ul{n}-T}K_i\right) \right] 
\end{align*}
   \item
If $\dgrm{X}$ is objectwise fibrant, there is a natural objectwise equivalence:
\begin{align*}
    \mathrm{spt}\left(\left(\bigwedge\limits_{i=1}^nR^{K_i}\right)_{\rm cof},\dgrm{X}\right)
                   &\simeq {\rm
                   hofib}\left[\dgrm{X}\left(\bigvee_{i=1}^nK_i\right)\to\holim_{T\in
                   \cal P_0(\ul{n})}\dgrm{X}\left(\bigvee_{\ul{n}-T}K_i\right) \right] \\ 
                   & \cong{\rm cr}_n\dgrm{X}(K_1,...\, ,K_n) 
\end{align*}
\end{punkt}
\end{lemma}

\begin{proof}
  Follows from the enriched Yoneda isomorphism $\mathrm{spt}(R^K,\dgrm{X})\cong\dgrm{X}(K)$
  in the same way as Lemma~\ref{cross effect}.
\end{proof}

The spectrum $\partial^{(n)}\dgrm{X}(*)$ for any homotopy functor $\dgrm{X}$ was introduced in 
\cite[p.~686]{Goo:calc3}; see also \cite[pp.~14-15]{Kuhn:overview}. There 
$\partial^{(n)}\dgrm{X}(*)$ is 
called the $n$-th derivative of $\dgrm{X}$ at $*$ and identified as ${\rm cr}_n\dgrm{X}(S^0,...\, ,S^0)$.
Permuting the zero-spheres induces a $\Sigma_n$-action on $\partial^{(n)}\dgrm{X}(*)$. To recover this
$\Sigma_n$-action on the left hand side in Lemma~\ref{lem:stable-cross-effect}, observe that
$\bigwedge_{i=1}^n R^K$ has a natural $\Sigma_n$-action
permuting the factors in the smash product. Choose a cofibrant replacement 
    $$ \id^n_{\rm cof} \rightarrow \bigwedge_{i=1}^n  R^{S^0} $$
in the naive $\Sigma_n$-equivariant model structure from Remark~\ref{naive}(ii). Then $\id^n_\mathrm{cof}$
is in particular a projectively cofibrant small functor which is weakly equivalent via a $\Sigma_n$-equivariant
map to the functor $K\mapsto K^{\wedge n}$. Taking $\Sigma_n$-orbits gives the canonical map
\[ K^{\wedge n}_{h\Sigma_n} \simeq \id^n_\mathrm{cof}(K)_{\Sigma_n} \rightarrow K^{\wedge n}_{\Sigma_n}\]
from homotopy orbits to orbits. Let $\rho_n\colon \Spec^{\cal{S}_*}\to \Spec^{\Sigma_n}$ be the functor
which maps $\dgrm{X}$ to the $\Sigma_n$-spectrum $\mathrm{spt}(\id^n_\mathrm{cof},\dgrm{X})$ whose action
is induced by the $\Sigma_n$-action on $\id^n_\mathrm{cof}$. 
From \ref{lem:stable-cross-effect} we deduce a natural $\Sigma_n$-equivariant weak equivalence 
    \begin{equation}\label{eq:4}
      \rho_n\dgrm{X}\simeq{\rm cr}_n\dgrm{X}(S^0,...\, ,S^0)\cong\partial^{(n)}\dgrm{X}(*).
    \end{equation}
   
The left adjoint $ \lambda_n\co\Spec^{\cal{S}_*}\to\Spec^{\cal{S}_*}$ of $\rho_n$ is given by
    $$ \lambda_nE := (\sm{E}{\id^n_{\rm cof}})_{\Sigma_n}. $$

\begin{proposition} \label{Quillen adjunction}
The functors $\lambda_n\co\Spec^{\Sigma_n}\leftrightarrows\Spec^{\cal{S}_*}:\!\rho_n$ form a Quillen pair,
where $\Spec^{\mathcal{S}_\ast}$ is equipped with the $n$-homogeneous model structure.
%
\end{proposition}

\begin{proof}
The functor $\rho_n$ maps objectwise trivial fibrations to trivial fibrations, since $\id^n_{\rm cof}$ is 
projectively cofibrant. Therefore $\lambda_n$ preserves cofibrations. If $E\to F$ is a stable equivalence of
cofibrant $\Sigma_n$-spectra, then $\sm{E}{\id^n_\mathrm{cof}} \to \sm{F}{\id^n_\mathrm{cof}}$ is an objectwise 
stable equivalence of objectwise cofibrant functors. The cofibrancy implies that taking $\Sigma_n$-quotients yields
an objectwise stable equivalence $\lambda_n(E) \to \lambda_n(F)$. Since $\Spec^{\Sigma_n}$ has as generating trivial cofibrations
a set of stable equivalences of cofibrant $\Sigma_n$-spectra, $\lambda_n$ is a left Quillen functor to the projective model
structure, hence also to the $n$-excisive model structure. 

To see that $\lambda_n$ is also a left Quillen functor to the $n$-homogeneous model structure, it suffices to check
that $\lambda_n$ maps any  cofibration $E\to F$ of $\Sigma_n$-spectra to an $n$-homogeneous cofibration. By the
dual of Corollary~\ref{cor:char-fib}, it remains to
prove that
\begin{equation}\label{eq:2}
\xymatrix{ M_n\lambda_n(E) \ar[r]\ar[d] & M_n\lambda_n(F) \ar[d] \\
        \lambda_n(E) \ar[r] & \lambda_n(F)}
\end{equation}
is a homotopy pushout square. The functor $K\mapsto \sm{E}{\id^n_\mathrm{cof}}$ is $n$-homogeneous by
\cite[Lemma 3.1]{Goo:calc3}. Taking orbits of spectra with free action is a homotopy colimit, thus
$\lambda_n$ takes values in $n$-homogeneous homotopy functors by~\ref{prop:dn-commutes}. Hence
the vertical arrows in Diagram~(\ref{eq:2}) are objectwise weak equivalences, which concludes the proof.
\end{proof}

\begin{remark}\label{rem:not-quillen}
  The existence of non-smashing localization functors implies that $\lambda_1$ is
  not a Quillen equivalence -- see \cite[Ex.~3.3]{Kuhn:overview}. 
  The problem is that any functor of the form $\lambda_n(E)$ is determined by its
  values on finite spaces. As one may deduce
  from \cite{Goo:calc3}, 
  the total left derived of $\lambda_n$ induces a full embedding of homotopy categories.
  In the last section, we will describe the image again as the homotopy category of
  a model category on small spectrum-valued functors.
\end{remark}

\section{Finitary functors}
\label{section:finitary}


Let $i\colon \cal{S}_*^{\mathrm{fin}} \hookrightarrow \mathcal{S}_\ast$ be the full inclusion of the category of finite pointed simplicial
sets. Enriched (over $\cal S_\ast$) left Kan extension along $i$ defines a functor
\[ i_\sharp \colon \cal{S_*}^{\cal{S}_*^{\rm fin}} \to \cal{S_*}^{\cal{S}_*} \]
having $i^\ast(\dgrm{X}) = \dgrm{X}\circ i$ as a right adjoint. Since $i$ is a full embedding,
the unit $\Id \to i^\ast \circ i_\sharp$ is a natural isomorphism.
Abbreviate $\mathbf{SF}\colon = \cal{S_*}^{\cal{S}_*^{\rm fin}}$.


In \cite{Lydakis}, Lydakis constructed three cofibrantly generated model structures
on $\mathbf{SF}$, the first (projective) model structure has levelwise
weak equivalences and fibrations. The second (homotopy functor) model
structure is obtained by localizing the projective model
structure with respect to those maps between representable
functors in $\mathbf{SF}$ that are induced by weak
equivalences of finite pointed simplicial sets. Fibrant objects in the
homotopy functor model structure are precisely the objectwise fibrant
homotopy functors. A fibrant replacement in the homotopy functor model
structure may be chosen as 

\[ \dgrm{X} \to \dgrm{X}^h \colon = \mathrm{fib}\circ i^\ast \bigl( (i_\sharp \dgrm{X}) \circ \mathrm{fib}\bigr).\] 

The third (stable) model structure coincides with the 1-excisive model structure, which we will construct soon.
Analogs of the projective and the homotopy functor model structures exist in the category 
$\mathbf{SpF}\colon = \Spec^{\mathcal{S}^\mathrm{fin}_\ast}$.
The homotopy functor model structure will be denoted by $\mathbf{SF}_h$ (resp.,~$\mathbf{SpF}_h$).

\begin{remark}
One could expect that the stable model structure should arise only
after a colocalization of the 1-excisive model structure into a
1-homogeneous model structure. However, Lydakis works in the category
of functors enriched over $\cal S_\ast$ and all functors in this
category are already reduced (cf., Remark~\ref{alternative}). As we
shall see bellow, the case $n=1$ is somewhat special; for bigger
values of $n$ the colocalization is necessary in order to model the
homotopy theory of $n$-homogeneous functors.
\end{remark}

\begin{proposition}\label{prop:lift}
  There is a cofibrantly generated proper
  model structure on the category $\cal{S_*}^{\cal{S}_*}$ (resp.,~$\Spec^{\mathcal{S}_\ast}$) such that
  $f$ is a weak equivalence (fibration) if and only if $i^\ast(f)$ is a weak equivalence (fibration) 
  in $\mathbf{SF}_h$ (resp.,~$\mathbf{SpF}_h$). The induced left Quillen functor $i_\sharp$ is a 
  Quillen equivalence.
\end{proposition}

\begin{proof}
  The existence of the cofibrantly generated model structure 
  follows from \cite[Theorem 2.1.19]{Hovey}. The conditions in \cite[Theorem 2.1.19]{Hovey} are easily checked, since
  the unit $\Id \to i^\ast i_\sharp$ is an isomorphism. The latter also implies the statement about the Quillen equivalence.
  Right properness is obvious, left properness follows from left properness of simplicial sets, 
  since any finitary cofibration is in particular an
  objectwise cofibration.
\end{proof}

The model structure from Proposition~\ref{prop:lift} will be called the \emph{finitary} homotopy functor model
structure. Its cofibrations are called \emph{finitary} cofibrations, fibrant replacement is still denoted $\dgrm{X}\to \dgrm{X}^h$. 
Recall from \cite[Def.~5.10]{Goo:calc3} that a homotopy functor is called \emph{finitary}
if it commutes with filtered homotopy colimits.
In the finitary homotopy functor model structure, every fibrant and cofibrant object is a finitary functor. Therefore, every small functor $\dgrm{X}$ is weakly 
equivalent to a finitary homotopy functor, namely the functor $i_\sharp (i^\ast \dgrm{X}^h)^\mathrm{cof}$.
To see this, it suffices to observe that any simplicial functor $\dgrm{X}$ defined on finite pointed simplicial sets
commutes with filtered colimits. This is true, since $\dgrm{X}$ is a colimit of representable functors $R^K$,
where $K$ is finite, and these representable functors commute with filtered colimits.

The endofunctor on the category of small functors $T_n^\infty$, used in the next lemma, was described in Diagram~(\ref{eq:1}). 
\begin{lemma}\label{lem:n-excisive-finite}
  The natural transformation $\dgrm{X} \to T^\infty_n (\dgrm{X}^h)$ satisfies axioms {\bf (A.4)}, {\bf (A.5)} and
  {\bf (A.6)} for the finitary homotopy functor model structure.  
\end{lemma}

\begin{proof}
  Suppose that $f\colon \dgrm{X}\to \dgrm{Y}$ is a map such that $f(K)$ is a weak equivalence
  for every finite space $K$. Then also the map $T_n(f)(K)$ is a weak equivalence for every finite space $K$,
  since $f(U\star K)$ is a weak equivalence for every finite set $U$. It follows that
  $T_n^\infty(f)(K)$ is a weak equivalence for every finite space $K$.
  This gives {\bf (A.4)}. Axioms {\bf (A.5)} and {\bf (A.6)} follow as 
  in the proof of Proposition~\ref{A4-6 again}. 
\end{proof}

Any map $f\colon \dgrm{X}\to \dgrm{Y}$ such that $i^\ast T^\infty_n(f^h)$ is an objectwise weak equivalence
will be called \emph{finitary} $n$-\emph{excisive} equivalence. 

\begin{proposition}\label{prop:n-excisive-finite}
  There is a finitary $n$-excisive model structure on $\mathcal{S}_\ast^{\mathcal{S}_\ast}$ resp.~$\Spec^{\mathcal{S}_\ast}$
  having finitary cofibrations as cofibrations and finitary $n$-excisive equivalences as weak equivalences.
  The identity functor $\Id$ is a left Quillen functor to the $n$-excisive model structure.
\end{proposition}

\begin{proof}
  The first statement follows from Lemma~\ref{lem:n-excisive-finite} and Theorem~\ref{thm:bousfield-machine}.
  For the second statement, observe that $i^\ast \circ T^\infty_n = T^\infty_n \circ i^\ast$. It follows that
  any $n$-excisive equivalence is in particular a finitary $n$-excisive equivalence. The analog
  of corollary~\ref{cor:char-fib} then shows that an $n$-excisive fibration is in particular a finitary
  $n$-excisive fibration, which completes the proof.
\end{proof}


A very similar model structure was constructed in \cite{Dwyer_localizations} on the category of functors from finite $CW$-complexes to topological spaces.

Let $M_n\dgrm{X} \to \dgrm{X}$ be the homotopy fiber of the natural map $\dgrm{X} \to T_{n-1}^\infty \dgrm{X}^h$.
Say that a map $f$ is a \emph{finitary} $n$-homogeneous equivalence if the map $M_n(f)$ is a finitary
$n$-excisive equivalence.

\begin{lemma}\label{lem:finitary-homogeneous}
  The natural transformation $M_n\to \Id_{\Spec^{\mathcal{S}_\ast}}$
  satisfies axioms {\bf (A.4)}${}^\op$, {\bf (A.5)}${}^\op$ and
  {\bf (A.6)}${}^\op$ for the finitary $n$-excisive model structure on ${\Spec^{\mathcal{S}_\ast}}$.
\end{lemma}

\begin{proof}
  Since the finitary $n$-excisive model structure is right proper, $M_n\dgrm{X}$ coincides
  with the homotopy limit in the finitary $n$-excisive model structure. Hence
  axiom {\bf (A.4)}${}^\op$ holds. For the remaining axioms consult the proof of Lemma~\ref{A5-dual} and
  Lemma~\ref{left properness for M_n}.
\end{proof}

\begin{proposition}\label{prop:finitary-homogeneous}
  There is a finitary $n$-homogeneous model structure on $\Spec^{\mathcal{S}_\ast}$
  having finitary $n$-homogeneous equivalences as weak equivalences and finitary $n$-excisive fibrations as fibrations.
  The identity functor $\Id$ is a left Quillen functor to the $n$-homogeneous model structure.
\end{proposition}  

\begin{proof}
  The existence of the finitary $n$-homogeneous model structure follows from Lemma~\ref{lem:finitary-homogeneous}
  and the dual of Theorem~\ref{thm:bousfield-machine}. It is then obvious from Proposition~\ref{prop:n-excisive-finite}
  that any $n$-homogeneous fibration
  is also a finitary $n$-homogeneous fibration. Suppose $f\colon \dgrm{X}\to \dgrm{Y}$ is a map such that
  $M_n(f)$ is an $n$-excisive equivalence. Then $M_n(f)$ is also a finitary $n$-excisive equivalence, by the
  proof of Proposition~\ref{prop:n-excisive-finite}. 
  Thus $f$ is a finitary $n$-homogeneous equivalence, and the result follows. 
\end{proof}

\begin{theorem}\label{thm:quillen-eq}
  The functor $\lambda_n\colon \Spec^{\Sigma_n} \to \Spec^{\mathcal{S}_\ast}$
  is a left Quillen equivalence to the finitary $n$-homogeneous model structure.
\end{theorem}

\begin{proof}
  The right adjoint $\rho_n$ of $\lambda_n$ is constructed as $\mathrm{spt}(\id^n_\mathrm{cof},-)$,
  where $\id^n_\mathrm{cof}$ is a $\Sigma_n$-equivariant cofibrant replacement of the
  functor $K\mapsto K^{\wedge n}$. Since this functor is a (finite) colimit of functors represented
  by finite spaces, the functor $\id^n_\mathrm{cof}$ is finitary cofibrant. As in the proof of Proposition~\ref{Quillen adjunction}
  one may conclude that $\lambda_n$ is a left Quillen functor.
  
  To show that $\lambda_n$ is actually a Quillen equivalence, it suffices to prove that 
  $\rho_n$ reflects weak equivalences of fibrant objects, and that the derived unit
  $E\to \rho_n(\lambda_n E)_\mathrm{fib}$ is a weak equivalence for every cofibrant
  $\Sigma_n$-spectrum $E$. The latter follows from the equivalence in equation~(\ref{eq:4}), 
  since $\lambda_nE$ is weakly equivalent to the functor $K\mapsto (E\wedge K^{\wedge n})_{h\Sigma_n}$.
  So let $f\colon \dgrm{X} \to \dgrm{Y}$ be a map
  of finitary $n$-excisively fibrant functors such that $\rho_n(f)$ is a weak equivalence
  of $\Sigma_n$-spectra. Because $\rho_n$ preserves weak equivalences of fibrant
  objects, we may assume that both $\dgrm{X}$ and $\dgrm{Y}$ are cofibrant in the finitary $n$-homogeneous
  model structure, thus in particular finitary and $n$-homogeneous. Via the equivalence displayed
  in equation~(\ref{eq:4}), we may deduce from
  \cite[Prop.~5.8 and Prop.~3.4]{Goo:calc3} that $f(K)$ is a weak equivalence for every
  finite space $K$, and thus a weak equivalence in the finitary $n$-homogeneous model structure.
\end{proof}

\begin{remark}
The statements analogous to Proposition~\ref{prop:n-excisive-finite},
Lemma~\ref{lem:finitary-homogeneous}, Proposition~\ref{prop:finitary-homogeneous},
and Theorem~\ref{thm:quillen-eq} are satisfied for the categories
$\mathbf{SF}_h$ and $\mathbf{SpF}_h$ instead of $\Spec^{\Sigma_n}$ and
$\Spec^{\mathcal{S}_\ast}$ respectively. The same methods of proof
work. This extends the results of Lydakis \cite{Lydakis} and gives an
alternative interpretation of Goodwillie's classification theorem. The
details are left to the interested reader.
\end{remark}

\bibliographystyle{abbrv}
\bibliography{Xbib}

\end{document}